\documentclass[letterpaper, 10 pt, conference]{ieeeconf}  

\IEEEoverridecommandlockouts                              

\usepackage{graphicx}
\usepackage{amsmath}
\usepackage{float}
\usepackage{amssymb}
\usepackage{url}
\usepackage{xcolor}
\usepackage[colorlinks, allcolors=blue]{hyperref}

\usepackage{algorithm}
\usepackage{algpseudocode}
\newtheorem{lemma}{Lemma}
\newtheorem{definition}{Definition}
\newtheorem{theorem}{Theorem}
\newtheorem{remark}{Remark}

\newcommand{\df}{\overset{\text{def}}{=}}
\newcommand{\defe}{\overset{\mathrm{def}}{=}}

\newcommand{\R}{\mathbb{R}} 
\newcommand{\N}{\mathbb{N}} 

\allowdisplaybreaks

\title{\LARGE \bf
Minimisation of Submodular Functions Using Gaussian Zeroth-Order Random Oracles
}

\author{Amir Ali Farzin$^{1}$, Yuen-Man Pun$^{1}$, Philipp Braun$^{1}$, Tyler Summers$^{2}$, and Iman Shames$^{3}$
\thanks{$^{1}$School of Engineering, Australian National University and CIICADA Lab
       {\tt\small \{amirali.farzin,philipp.braun, yuenman.pun\}@anu.edu.au}}%
\thanks{$^{2}$ Department of Mechanical Engineering, University of Texas at Dallas
       {\tt\small tyler.summers@utdallas.edu}}%
\thanks{$^{3}$ Department of Electrical and Electronic Engineering, University of Melbourne and CIICADA Lab
       {\tt\small iman.shames@unimelb.edu.au}}%
}

\begin{document}

\maketitle
\thispagestyle{empty}
\pagestyle{empty}

\begin{abstract}
We consider the minimisation problem of submodular functions and investigate the application of a zeroth-order method to this problem. The method is based on exploiting a Gaussian smoothing random oracle to estimate the smoothed function gradient. We prove the convergence of the algorithm to a global $\epsilon$-approximate solution in the offline case and show that the algorithm is Hannan-consistent in the online case with respect to static regret. Moreover, we show that the algorithm achieves $O(\sqrt{NP_N^\ast})$ dynamic regret, where $N$ is the number of iterations and $P_N^\ast$ is the path length. The complexity analysis and hyperparameter selection are presented for all the cases. The theoretical results are illustrated via numerical examples.
\end{abstract}

\section{INTRODUCTION}

In this paper, we study a minimisation problem of the form
	\begin{equation}\label{mainp}
		\underset{S \subseteq Q}{\min}\;f(S),
	\end{equation}
	where the objective function $f\colon 2^{Q}\rightarrow\mathbb{R}$ is a set function over a ground set $Q$ of $n$ elements, where $2^{Q}$ denotes the power set of $Q$. Without loss of generality, we assume $Q = [n] = \{1,2,\dots,n\}$. We assume a minimiser for $f$ exists (which may not be unique) and consider computing an $\epsilon$-approximate minimiser of $f,$ i.e., $S\subseteq Q$ with $f(S)\leq f^\ast+\epsilon,$ where $f^\ast\df \min_{S \subseteq Q}\;f(S).$
 More specifically, 
 we focus on submodular function minimisation (SFM), which counts as a foundational problem in combinatorial optimisation. After the seminal work of \cite{edmonds1970submodular}, SFM found its place as a key tool in many areas. Its application can be seen in machine learning \cite{krause2010submodular}, economics \cite{topkis1998supermodularity}, computer vision \cite{jegelka2013reflection}, and speech recognition \cite{lin2011optimal}. For more 
 information, we refer to the survey papers
 \cite{iwata2008submodular, mccormick2005submodular}. In the seminal work 
 \cite{grotschel1981ellipsoid}, the authors showed that SFM problems
 can be solved 
 in polynomial time. Since then, there has been extensive research and multiple advances in characterising SFM \cite{chakrabarty2017subquadratic,chakrabarty2014provable,axelrod2020near}.

 An important key in solving SFM is the fact that SFM reduces to minimising the Lovász extension \cite{lovasz1983submodular} over $[0,1]^n.$ 
 The Lovász extension enjoys desirable properties such as convexity and Lipschitz continuity, which can be leveraged to solve SFM. 
 The Lovász extension is non-smooth and thus 
 gradient methods can not be used to find its minimisers. Hence, leveraging the subgradient of the Lovász extension plays a key part in the developed of methods to solve both offline SFM \cite{hamoudi2019quantum,axelrod2020near} and online SFM \cite{hazan2012online,jegelka2011online,cardoso2019differentially}. Since the Lovász extension is non-smooth, zeroth-order methods are a suitable choice to solve SFM.

Optimisation frameworks that rely only on the evaluations of the objective function (but not on its derivatives) are commonly called zeroth-order (ZO) or derivative-free methods \cite{audet2017introduction}. Multiple ZO algorithms have been developed to solve various optimisation problems. The authors in \cite{nesterov_random_2017}, proposed and analysed a method based on ZO oracles leveraging Gaussian smoothing. This method has been extensively studied for different classes of functions in both offline \cite{farzin2024minimisation,farzin2025minimisation} and online settings \cite{shames2019online}. In this work, we leverage a Gaussian smoothing framework to solve SFM in both offline and online settings. Leveraging this method, we do not need to calculate the subgradient of the Lovász extension.
In this work, we first consider the offline SFM and show that the Gaussian smoothing ZO method, in the expectation sense, converges to an $\epsilon$-approximate solution 
in $O(n^2\epsilon^{-2})$ function calls. Next, we study 
online SFM under both static and dynamic regrets (unlike the previous works on online SFM, which only consider the static regret). We show that, in the expectation sense, the algorithm is Hannan-consistent (see \cite{hazan2012online}) in terms of static regret, meaning that its average regret per iteration vanishes as the number of iterations grows. Moreover, it achieves a dynamic regret of $O(\sqrt{NP_N^\ast})$, where $N$ is the number of iterations and $P_N^\ast$ denotes the path length of the minimisers. It is known that the dynamic regret bound is dependent on the path length \cite{zinkevich2003online,hazan2016introduction}.

\emph{Outline:} The paper is organised as follows. Preliminaries are introduced in Section~\ref{sec:prem}. In Section~\ref{sec:main}, the main framework is given, and convergence and complexity results related to the proposed framework are presented
for offline and online settings. 
Illustrative examples are given in Section~\ref{sec:exm}. Lastly, we conclude our paper and discuss potential future directions in Section~\ref{sec:conc}. 

\emph{Notation:}
For $x,y\in \mathbb{R}^d$,  $\langle x , y \rangle=x^\top y$ and let $\|\cdot\|$ be the Euclidean norm of its argument. The gradient of a differentiable function $f: \mathbb{R}^{n+m} \to \mathbb{R}$ is denoted by  $\nabla f$. Let $\{0,1\}^n \subseteq \N^n$
denote the set of $n$-dimensional binary vectors, where each coordinate is either $0$ or $1$.  There is a natural bijection between a vector $x \in \{0,1\}^n$ and a subset $S \subseteq [n]$. For a set $S \subseteq [n]$, let $\chi_S \in \{0,1\}^n$ denote its characteristic vector, 
defined by $\chi_S(i) = 1$ if $i \in S$ and $\chi_S(i) = 0$ otherwise.
 We will freely move between these two representations and, in particular, interchangeably use
\begin{align}
    f(S), \ S \subseteq [n] \quad \text{and} \quad f(\chi_S), \ \chi_S \in  \{0,1\}^n.
\end{align}
The function $F:D\rightarrow \mathbb{R}$ with $D\subseteq \mathbb{R}^n$ is Lipschitz if there exists a \emph{Lipschitz constant} $L_0>0$, such that
	\begin{equation}\label{def:lip}
		|F(x)-F(y)|\leq L_0\|x-y\|,\quad \forall \ x,y\in D. 
	\end{equation}
\section{Preliminaries}\label{sec:prem}
In this section, we review the necessary definitions and
background material of submodular functions, Lovász extension, and Gaussian smoothing.

\subsection{Submodular Functions and the Lovász Extension}

In this section, we give the background related to submodular functions and their  Lovász Extensions.

\begin{definition}[Submodular Functions \cite{hazan2012online}]\label{def:sub}
    A function $f: 2^{[n]} \to \mathbb{R}$ is called \emph{submodular} if it exhibits diminishing marginal returns, i.e., 
for all $S \subseteq T \subseteq [n]$ and any 
$i \notin T$,  
\[
    f(S \cup \{i\}) - f(S) \;\;\ge\;\; f(T \cup \{i\}) - f(T).
\]  
An equivalent characterisation is that for all $S, T \subseteq [n]$,  
\[
    f(S) + f(T) \;\;\ge\;\; f(S \cap T) + f(S \cup T).
\]  
\end{definition}
    
Let $f : 2^{[n]} \to \mathbb{R}$ be a submodular function defined on subsets of a ground set $E$.  
We assume that access to $f$ is provided only through a value oracle. That is, given any set $S \subseteq [n]$, the oracle returns the value $f(S)$.
The 
\emph{Lovász extension} provides a continuous, convex extension of a submodular function 
$f : 2^{[n]} \to \mathbb{R}$ to a function $f^L: 
[0,1]^n \to \mathbb{R}$, where $ [0,1]^n\subset \R^n$ is the unit hypercube.
\begin{definition}[Lovász Extension \cite{hazan2012online}]\label{def:lovae}
    Let $x\in[0,1]^n$ and
    let $A_i\subset [n]$ ($i\in [n]$, $n\in \N$), be a chain of subsets
    $A_0\subset A_1 \subset \cdots \subset A_p$ 
    such that $x$ can be expressed as a convex combination $x=\sum_{i=0}^n \lambda_i\chi_{A_i}$ where $\lambda_i>0$ and $\sum_{i=0}^n\lambda_i=1.$ Then, 
    the Lovász extension $f^L$ at $x$ is defined to be
    \begin{align}\label{eq:eqle1}
        f^L(x) =  \sum_{i=0}^p \lambda_if(A_i).
    \end{align}
\end{definition}
\begin{remark}  \label{rem:lovae}
According to \cite{hazan2012online}, the Lovász extension can alternatively be obtained by
sampling 
a threshold $\tau\in[0,1]$ uniformly at random, and considering 
    \begin{align}\label{eq:eqle2}
        f^L(x) = E_\tau[f(S_\tau)], \qquad S_\tau=\{i\in [n]:x(i)>\tau\}.
    \end{align}
\end{remark}

Remark~\ref{rem:lovae}, 
implies that for all sets $S\subseteq[n],$ it follows that $f^L(\chi_S)=f(S).$
We use the following 
properties of the Lovász extension function shown in \cite{lovasz1983submodular,fujishige2005submodular,bach2013learning,ando2002k}.
\begin{lemma}\label{lm:lmlova}
    Let $f:2^{[n]}\rightarrow \mathbb{R}$ be a submodular function and $f^L$ be its Lovász extension. Then, we have that (i) $f^L$ is convex, (ii) $f^L$ is piecewise linear,  
    (iii) $\min_{S\subseteq[n]} f(x) = \min_{x\in\{0,1\}^n} f^L(x)=\min_{x\in[0,1]^n} f^L(x)$, and (iv) the set of minimisers of $f^L(x)$ in $[0,1]^n$ is the convex hull of the set of minimisers of $f^L(x)$ in $\{0,1\}^n$. 
\end{lemma}
Also, in \cite[Lemma 2.2]{axelrod2020near}, the authors show
that we can go from an approximate minimiser of $f^L$ to an approximate minimiser of $f$ in $O(n)$ function calls.  Thus, leveraging this fact and considering Lemma~\ref{lm:lmlova}, we can consider solving the non-smooth convex minimisation problem of the form
\begin{equation}\label{mainp2}
		\underset{x \in [0,1]^n}{\min}\;f^L(x),
\end{equation}
instead of Problem~\ref{mainp}. Since $f^L$ is continuous and the domain is compact and nonempty, at least a minimiser for $f^L$ exists, which may not be unique. Since
$f^L$ is non-smooth, ZO methods are good candidates to solve \eqref{mainp2}. 

\subsection{Gaussian smoothing}
In this section, following \cite{nesterov_random_2017}, we define the Gaussian smoothed version of a continuous function $F:\mathbb{R}^n\rightarrow\mathbb{R},$ termed $F_\mu$ as below:
\begin{align}\label{eq:gsmooth}
	\begin{split}
		F_\mu(x) &\defe \frac{1}{\kappa}\underset{E}{\int}F(x+\mu u)e^{-\frac{1}{2}||u||^2}\mathrm{d}u,\\
        \kappa &\defe\underset{E}{\int}e^{-\frac{1}{2}||u||^2}\mathrm{d}u=\frac{(2\pi)^{n/2}}{[\det B]^{\frac{1}{2}}},
	\end{split}
\end{align}
where vector $u \in \R^n$ is sampled from $\mathcal{N}(0,B^{-1})$, $B\in \R^{n\times n}$ is positive definite, and $\mu \in \R_{>0}$ is the 
smoothing parameter. In \cite{nesterov_random_2017}, it is shown that whether $F$ is differentiable or not, $F_\mu$ is always differentiable, and
$\nabla F_\mu(x) = E_u[\frac{F(x+\mu u)}{\mu}u].$ 
Instead of this one-point gradient estimation, to reduce the variance of the estimator, we consider the two-point estimator and  define the random oracle $g_\mu$ as
\begin{equation}\label{eq:eq23}
	g_\mu(x) \defe \tfrac{1}{\mu} (F(x+\mu u)-F(x) Bu),
\end{equation}
where $u$ and $B$ are defined above. 

\section{Main Results}\label{sec:main}
In this section, we investigate the application of a ZO algorithm in solving offline and online minimisation of submodular functions by considering their Lovász extensions. First, we show that we can use Gaussian smoothing on the Lovász extension $f^L(x).$ Note that it is not trivial to apply~\eqref{eq:eq23} with $f^L$  
as 
$x+\mu u$
can be outside of the unit hypercube. To use Gaussian smoothing, we need to extend
the domain of $f^L$ to $\R^n$ using \cite[Definition 3.1]{bach2013learning}.
\begin{definition}[Extended domain Lovász Extension]\label{def:lovaee}~Giv{\-}en a set-function $f:2^{[n]} \rightarrow\R$ such that $f(\emptyset) = 0$, the Lovász extension $f^L: \mathbb{R}^n \rightarrow \mathbb{R}$ is defined as follows:
    For $x \in \mathbb{R}^n$, order the components in decreasing order 
    $x({j_1}) \geq \cdots \geq x({j_n})$, where $(j_1, \ldots, j_n)$
    is a permutation. The Lov\'{a}sz extension is given by
\begin{align}\label{eq:extlov}
f^L(x) \hspace{-0.05cm} = \hspace{-0.05cm} \sum_{k=1}^n x_{j_k} \hspace{-0.05cm} \left[ f(\{j_1, \ldots, j_k\}) - f(\{j_1, \ldots, j_{k-1}\}) \right]
\end{align}
\end{definition}
\begin{remark}
According to \cite[(3.3)]{bach2013learning}, it holds that \eqref{eq:eqle2} and \eqref{eq:extlov} coincide on the domain $[0,1]^n$ and thus the notation $f^L$ can be used in Definition \ref{def:lovae} and in Definition \ref{def:lovaee}. From here, we use $f^L$ to denote \eqref{eq:extlov}.
\end{remark}
 Letting $f^L_\mu=E_{u\sim\mathcal{N}(0,I)}[f^L(x+\mu u)],$ 
 the two-point gradient estimator given in~\eqref{eq:eq23} with $F=f^L$ gives an unbiased estimation of $\nabla f^L_\mu(x)$ for solving~\eqref{mainp2}.
 In the next section, we will analyse the convergence for the offline problem~\eqref{mainp2} and extend the results to the online case.

\subsection{Offline Problem}

In this section, first we introduce Algorithm~\ref{alg:ZOGD}. This framework has been studied for different minimisation problems in the literature \cite{nesterov_random_2017,farzin2024minimisation,farzin2025minimisation}.
\begin{algorithm}[htb]
	\caption{ZO-GD}\label{alg:ZOGD}
	\begin{algorithmic}[1]
		\State Input: $x_0\in\mathcal{K},\;\{h_{k}\}_{k=0}^N\subset \mathbb{R}_{>0},\;\mu>0,\;N,t_k\in\mathbb{N}$
		\For {$k = 0,\dots, N$}
		\State Sample $u_k^0, \dots, u_k^{t_k}$ from $\mathcal{N}(0, \mathbb{I})$
            \State Obtain $g_{\mu}^0(x_k), \dots, g_{\mu}^{t_k}(x_k)$ using $u_k^0, \dots, u_k^t$ and \eqref{eq:eq23}
            \State $g_{\mu}(x_k) = \frac{1}{t_k} \sum_{i=0}^{t_k} g_{\mu}^i(x_k)$ \label{line:6}
		\State $\bar{x}_{k+1} = x_k - h_kg_\mu(x_k)$ \label{line:5}
		\State $x_{k+1} = \mathrm{Proj}_{[0,1]^n}(\bar{x}_{k+1})$
		\EndFor
		\State return $x_{k}$ for all $k = 0,\dots, N.$
	\end{algorithmic}
\end{algorithm}
In each iteration of Algorithm~\ref{alg:ZOGD}, first,
to update our guess, we sample $t_k$ number of directions from the standard Gaussian distribution and calculate the random oracles, with $F=f^L,$ using the sampled directions and leverage the average of the oracles. It is easy to see that with $t_k$ number of samples instead of one, still $E[g_\mu(x_k)] = \nabla f^L_\mu(x_k),$ while the variance of $g_\mu(x_k)$ will be reduced \cite{farzin2025minmax}. Then, Algorithm \ref{alg:ZOGD} performs
a descent step using the calculated random oracle with $h_k$ as the step size. At last, we project back to the unit hypercube. The projection operator $\mathrm{Proj}_{[0,1]^n}:\R^n \rightarrow [0,1]^n$
is defined by 
\[
\mathrm{Proj}_{[0,1]^n,i}(x) =  \max\{0,\min\{1,x(i)\}\}, \quad i\in [n].
\]
Considering Lemma~\ref{lm:lmlova}, it can be concluded that $f^L$ is continuous and piecewise affine. Thus, $f^L$ is Lipschitz continuous with Lipschitz constant $L_0\in \R_{>0}$ defined through the maximal slope of the piecewise affine functions.
Following \cite{hazan2012online}, if in each iteration of Algorithm~\ref{alg:ZOGD} we
choose a threshold $\tau\in[0,1]$ uniformly at random, we can define the set 
\begin{equation}\label{eq:Sk}
    S_k=\{i\in[n]:x_k(i)>\tau\}.
\end{equation}
Next, we present the main theorem of this section.

\begin{theorem}\label{th:off}
   Consider Algorithm \ref{alg:ZOGD} and let $f:2^{[n]}\to\mathbb{R}$ be a submodular function and $f^L$ be its Lovász extension with $x^\ast$ as a
   minimiser. Let $L_0$ be the Lipschitz constant of $f^L$, $N\geq0$ be the number of iterations, $t_k=t\geq1$ be the number of samples, $r_0=\|x_0-x^\ast\|$, $\mu>0$ be smoothing parameter, $\mathcal{U}_k = [u_0,u_1,,\cdots,u_k]$, 
     $k\in [N]$.
     Moreover, let $\{x_k\}_{k\geq0}$ be the sequences generated by Algorithm~\ref{alg:ZOGD},     
     Then, for any iteration $N$, with \(h_{k}= h>0\),
     we have
	\begin{align}
    \begin{split}
		\frac{1}{N+1} \sum_{k=0}^{N}&E_{\mathcal{U}_k}[f^L(x_k)]-f^L(x^\ast) \\ &\leq\frac{n}{2h(N+1)}+\frac{\mu}{2} L_0n^{\frac{1}{2}}L_0^2h(n+4)^2.
        \end{split}\label{eq:eqoff}
	\end{align}
    Moreover, let $\epsilon>0$, $\mu\leq\frac{\epsilon}{2L_0n^{\frac{1}{2}}}$ and
    \begin{align*}
      h = \frac{r_0}{(N+1)^{\frac{1}{2}}L_0(n+4)},   \quad N\geq\Big\lceil \frac{4r_0^2L_0^2(n+4)^2}{\epsilon^2}-1\Big\rceil.
    \end{align*}
Then, we have 
    \begin{align}\label{eq:eqoff2}
		E_{\mathcal{U}_{N-1}}[f^L(\hat{x}_N)]-f^L(x^\ast)\leq\epsilon,
	\end{align}
    where $\hat{x}_N \overset{\mathrm{def}}{=} \arg\underset{x}{\min}[f^L(x):x\in\{x_0,\dots,x_N\}]$, and
     \begin{align}\label{eq:corooff}
		E_{\tau,\mathcal{U}_{N-1}}[f(\hat{S}_N)]-f(S^\ast)\leq\epsilon,
	\end{align}
    where $S^\ast$ is a minimiser of $f,$ $\tau$ is the threshold used in~\eqref{eq:Sk} and 
    $
    \hat{S}_N \overset{\mathrm{def}}{=} \underset{S\in\{S_0,\dots,S_N\}} {\arg\min} f(S).
    $
    
\end{theorem}

\begin{proof}
    The proof follows the proof of \cite[Thm. 6]{nesterov_random_2017}. Let $r_k = \|x_k-x^\ast\|.$ Using the non-expansiveness property of the projection to convex sets, we have
    \begin{align}\label{eq:th1eq1}
        r_{k+1}^2\leq&\|\bar{x}_{k+1}-x^\ast\|^2\leq\|x_k-x^\ast-hg_\mu(x_k)\|^2\notag\\\leq&r_k^2-2h\langle g_\mu(x_k),x_k-x^\ast\rangle+h^2\|g_\mu(x_k)\|^2.
    \end{align}
    Taking the expectation of 
    \eqref{eq:th1eq1} with respect to $u_k$ and considering the fact that similar to the proof of \cite[Thm. 4.1]{nesterov_random_2017}, we can obtain an upper-bound for \eqref{eq:eq23}, we have
    \begin{align*}
        E_{u_k}[r_{k+1}^2]\leq&r_k^2-2h\langle \nabla f^L_\mu(x_k),x_k-x^\ast\rangle+h^2L_0^2(n+4)^2.
    \end{align*}
Considering the convexity of $f^L$ (which implies that 
$f_\mu^L$ is convex \cite{nesterov_random_2017}) 
and \cite[Thm. 1]{nesterov_random_2017}, we have
    \begin{align}\label{eq:th1eq2}
        f^L&(x_k)-f^L(x^\ast) \leq\tfrac{r_k^2-E_{u_k}[r_{k+1}^2]}{2h}+\mu L_0n^{\tfrac{1}{2}}+\tfrac{hL_0^2(n+4)^2}{2}.
    \end{align}
    Considering the fact that $r_0^2=\|x_0-x^\ast\|^2\leq n,$ taking expectation with respect to $\mathcal{U}_{k-1}$ from \eqref{eq:th1eq2}, summing from $k=0$ to $k=N,$ and dividing by $N+1,$ we arrive at \eqref{eq:eqoff}. 
    If we substitute the 
    values in Theorem \ref{th:off} for $\mu,$ $h,$ and $N$ in \eqref{eq:eqoff}, we have
\begin{align*}
		\frac{1}{N+1} \sum_{k=0}^{N}E_{\mathcal{U}_k}[f^L(x_k)]-f^L(x^\ast)\leq&\; \epsilon.
	\end{align*}

We know that for $k\in [N]$, $f^L(x_k)-f(x^\ast)\geq0$. 
Considering the fact that if the average of a positive sequence is less than or equal to $\epsilon$, then the minimum is less than or equal to $\epsilon,$ i.e., we arrive at \eqref{eq:eqoff2}.    

Considering Lemma~\ref{lm:lmlova}, we know $f(S^\ast) = f^L(x^\ast)$ where $x^\ast$ is a
    minimiser of the Lovász extension. From
    the projection step in Algorithm~\ref{alg:ZOGD}, we know that $x_k\in[0,1]^n$ for all $k \in [N].$ Thus, considering Definition~\ref{def:lovae}, we have $f^L(x_k) = E_\tau[f(S_k)].$ Moreover, we know that $\tau$ and $u_k$ are independent random variables.
    Considering these facts and taking expectation from $f(S_k)-f(S^\ast)$ with respect to $\tau,$ then taking expectation with respect to $\mathcal{U}_{k-1}$, summing from $k=0$ to $k=N,$ and dividing by $N+1,$ and letting the hypothesis of Theorem~\ref{th:off} be satisfied, we have 
    \begin{align}
        \frac{1}{N+1} \sum_{k=0}^{N}E_{\tau,\mathcal{U}_k}[f(S_k)]-f(S^\ast)\leq\epsilon.
    \end{align}
    Similarly, we can pick the best guess and arrive at \eqref{eq:corooff}.
\end{proof}

Theorem~\ref{th:off}, shows that in the expectation sense, there exists an $\epsilon$-optimal point for $f^L$ in the sequence generated by Algorithm~\ref{alg:ZOGD}.

 \begin{remark}
     If instead of \eqref{eq:eq23}, we use the centeral difference or backward random oracle as $\hat{g}_\mu(x) = \frac{F(x+\mu u)-F(x-\mu u)}{2\mu}u$ or $\bar{g}_\mu(x) = \frac{F(x)-F(x-\mu u)}{\mu}u,$ we can obtain the same bounds
     as in Theorem~\ref{th:off}. This is because the upper bounds on the expectation of the square norm of all these oracles are the same for a Lipschitz continuous function.
 \end{remark}

\subsection{Online Problem}
In this section, we analyse the online minimisation of submodular functions using a ZO method. 
In the online submodular minimisation problem, over a sequence of iterations $k\in \N$,, 
an online decision maker repeatedly selects a subset \( S_k \subseteq [n] \). After choosing \( S_k \), the cost is determined by a submodular function \( f_k : 2^{[n]} \rightarrow \mathbb{R} \), and the decision maker incurs the cost \( f_k(S_k) \). We will analyse both the static regret and dynamic regret of Algorithm~\ref{alg:ZOGD}. 
\begin{definition}[Static Regret]
The static regret of the decision maker is defined as:
\[
\text{Regret}_N := \sum_{k=0}^N f_k(S_k) - \min_{S \subseteq [n]} \sum_{k=0}^N f_k(S).
\]
\end{definition}

If the sets \( S_k \) are chosen by a randomised algorithm, the expected regret over the randomness in the algorithm is considered. An online algorithm to choose the sets $S_k$ is said to be Hannan-consistent if it ensures that \(E[\text{Regret}_N] \leq o(N)\).
Here, as considered in \cite{shames2019online}, we assume that the submodular cost function can change between two function evaluations that are needed to calculate \eqref{eq:eq23}. Thus we consider working with a family of submodular functions $f_k: 2^{[n]} \rightarrow \mathbb{R}$ where $k\in\mathbb{N}\cup\{j+\frac{1}{2}|j\in\mathbb{N}\} = \{0,1/2,1,3/2,\dots\}.$  For simplicity, we denote $k+\frac{1}{2}$ by $k^+$ for all $k\in\mathbb{N}.$ Also, we denote the Lovász extension of $f_k$ by $f^L_k.$ We know that each function $f^L_k$ is Lipschitz with $L_{0,k}$ as the Lipschitz constant. Next,
We need to alter the random oracle defined in \eqref{eq:eq23} to adapt it to this setting. As mentioned after \eqref{eq:gsmooth} and before \eqref{eq:eq23}, 
we know that $\nabla f^L_{\mu,k}(x) = E_u[\frac{f^L_{k}(x+\mu u)}{\mu}u]$ and $E[u]=0,$ thus we can define 
\begin{align}\label{eq:gon}
    g_{\mu,k^+}(x) = \tfrac{1}{\mu} (f^L_{k}(x+\mu u)-f^L_{k^+}(x)) u
\end{align}
with
$\nabla f^L_{\mu,k}(x) = E_u[g_{\mu,k^+}(x)].$ Thus, to solve this problem, we consider Algorithm~\ref{alg:ZOGD} with below update step
\begin{align}\label{eq:newup}
    \bar{x}_{k+1} = x_k -h_kg_{\mu,k^+}(x_k).
\end{align}
\begin{theorem}\label{th:on}
For $k\in\{0,\ldots,N\}$, let $f_k:2^{[n]}\to\mathbb{R}$ be submodular functions and $f^L_k$ be their corresponding Lovász extension with $x^\ast$ as a minimiser of $\sum_{k=0}^N f^L_k(x)$. Let $\bar{L}_0 = \max_{k}L_{0,k}$ be the Lipschitz constant, $N\geq0$ be the number of iterations, $t_k=t\geq1$ be the number of samples, $\mu>0$ be the  smoothing parameter in \eqref{eq:gsmooth} and $\mathcal{U}_k = [u_0,u_1,,\cdots,u_k]$. 
     Moreover, let $\{x_k\}_{k\geq0}$ be the sequences generated by Algorithm~\ref{alg:ZOGD} with \eqref{eq:gon} and \eqref{eq:newup} as the update step in Lines \ref{line:6} and \ref{line:5}.    
     Then, for any iteration $N$, with \(h_{k}= h = \frac{n^{\frac{1}{2}}} {(N+1)^{\frac{1}{2}}L_0(n+4)}\) and \(\mu\leq\frac{1}{L_0n^{\frac{1}{2}}(N+1)^{\frac{1}{2}}}\),
     we have
	\begin{align}
    \begin{split}
	\sum_{k=0}^{N}E_{\mathcal{U}_k}[f_k^L(x_k)]&-\min_{x\in\mathcal{K}}\sum_{k=0}^{N}f_k^L(x)\\
    &\leq(N+1)^{\frac{1}{2}}\Big(1+n^{\frac{1}{2}}(n+4)\bar{L}_0\Big).
	\end{split} \label{eq:eqon} 
    \end{align}
    Moreover, we get 
    \begin{align}\label{eq:eqon2}
		E_{\tau,\mathcal{U}_{N-1}}[\text{Regret}_N]= O(n^{\frac{3}{2}}\sqrt{N}),
    \end{align}
    where $\tau$ is the threshold defined in Algorithm~\ref{alg:ZOGD}
\end{theorem}
\begin{proof}
    Similar to the steps of the proof of Theorem~\ref{th:off}, letting $r_k=\|x_k-x^\ast\|,$ we have 
    \begin{align}\label{eq:th2eq1}
        E_{u_k}[r_{k+1}^2]\leq&r_k^2-2h\langle \nabla f^L_{\mu,k}(x_k),x_k-x^\ast\rangle+h^2\bar{L}_0^2(n+4)^2
    \end{align}
Considering the convexity of $f^L_k$, which implies that 
$f_{\mu,k}^L$ is convex,
and \cite[Thm. 1]{nesterov_random_2017}, we have
    \begin{align}\label{eq:th2eq2}
        f_k^L&(x_k)-f_k^L(x^\ast) \notag\\
        &\leq \frac{r_k^2-E_{u_k}[r_{k+1}^2]}{2h}+\mu \bar{L}_0n^{\frac{1}{2}} +\frac{h\bar{L}_0^2(n+4)^2}{2}.
    \end{align}
    Considering the fact that $\|x_0-x^\ast\|^2\leq n,$ taking expectation with respect to $\mathcal{U}_{k-1}$ from \eqref{eq:th2eq2}, summing from $k=0$ to $k=N,$ and substituting the given values for $\mu,$ $h,$ and $N$, we arrive at \eqref{eq:eqon}.
    
    Considering Lemma~\ref{lm:lmlova}, $\sum_{k=0}^N f_k(S^\ast) = \sum_{k=0}^N f_k^L(x^\ast)$ where $S^\ast$ is the minimiser of $\sum_{k=0}^N f_k(S).$ Also, considering the projection step in Algorithm~\ref{alg:ZOGD}, we know that $x_k\in[0,1]^n$ for all $k \in \{0,\ldots,N\}.$ Thus, considering Definition~\ref{def:lovae}, we have $f_k^L(x_k) = E_\tau[f_k(S_k)].$ Moreover, we know that $\tau$ and $u_k$ are independent random variables.
    Considering these facts and taking expectation from $f_k(S_k)-f_k(S^\ast)$ with respect to $\tau,$ taking expectation with respect to $\mathcal{U}_{k-1}$, summing from $k=0$ to $k=N,$ and considering \eqref{eq:eqon}, we arrive at \eqref{eq:eqon2}.
\end{proof}

Theorem~\ref{th:on} implies that Algorithm~\ref{alg:ZOGD} with \eqref{eq:gon} and \eqref{eq:newup} as the update step in Lines \ref{line:6} and \ref{line:5}, respectively, is Hannan-consistent.
We should note that if $f_k = f_{k^+},$ we can arrive at the same results as in Theorem~\ref{th:on}. This is because we would still have $\nabla f^L_{\mu,k}(x) = E_u[g_{\mu,k^+}(x)].$ In the next remark, we discuss
what 
happens if, the random oracle is defined as 
\begin{align}\label{eq:gon2}
    \tilde{g}_{\mu,k^+}(x) = \frac{f^L_{k^+}(x+\mu u)-f^L_{k}(x)}{\mu}u
\end{align}
instead of \eqref{eq:gon}.
\begin{remark}\label{rem:reverse}
    If in the update step for the online setting,  
    \eqref{eq:gon2} is used instead of \eqref{eq:gon}, then we get $\nabla f^L_{\mu,k^+}(x) = E_u[\tilde{g}_{\mu,k^+}(x)].$ In this case, $\nabla f^L_{\mu,k^+}(x)$ is not an unbiased estimator
    of $\nabla f^L_{\mu,k}(x)$ and an extra assumption is necessary. In particular, we assume that 
    $|f^L_k(x)-f^L_{k^+}(x)|\leq V$ for all $x$ and where $V\in \R_{\geq 0}$ denotes a constant. This means that the change of the functions' values during the sampling period in each iteration is upper-bounded. With this assumption, we get the estimate
\begin{align}\label{eq:eqonrem}
	\sum_{k=0}^{N}&E_{\mathcal{U}_k}[f_k^L(x_k)]-\min_{x\in\mathcal{K}}\sum_{k=0}^{N}f_k^L(x)\\
    &\leq(N+1)^{\frac{1}{2}}\Big(1+\notag n^{\frac{1}{2}}(n+4)\bar{L}_0\Big)+2V(N+1),
	\end{align}
    and
    \(
		E_{\tau,\mathcal{U}_{N-1}}[\text{Regret}_N]= A +B,
    \)
where $A = O(n^{\frac{3}{2}}\sqrt{N})$ and $B = O(NV).$
\end{remark}

Next, we consider tracking the minimum value of each of $f_k.$ To this end, we
need to define the dynamic regret.

\begin{definition}[Dynamic Regret]
The dynamic regret of the decision maker is defined as:
\[
\text{D-Regret}_N := \sum_{k=0}^N f_k(S_k) - \sum_{k=0}^N \min_{S \subseteq [n]} f_k(S).
\]
\end{definition}
If the sets $ S_k\subseteq 2^{[n]}$ are chosen by a randomised algorithm, the expected regret over the randomness in the algorithm is considered. For the analysis of the dynamic regret,
following \cite{zinkevich2003online,hazan2016introduction}, we need to define the total path length 
\begin{align}\label{eq:pathl}
    P_N(u_0,\dots,u_N)=\sum_{i=1}^N\|u_i-u_{i-1}\|,
\end{align}
which allows us to
present our main result for the dynamic regret, corresponding to the
sequence generated by Algorithm~\ref{alg:ZOGD} with \eqref{eq:gon} and \eqref{eq:newup} as the update steps in Lines \ref{line:6} and \ref{line:5}.
\begin{theorem}\label{th:ond}
For $k\in \{0,\ldots,N\}$,
let $f_k:2^{[n]}\to\mathbb{R}$ be submodular functions and $f^L_k$ be their corresponding Lovász extension with minimiser $x^\ast_k$.
Let $\bar{L}_0 = \max_{k}L_{0,k}$ be the common Lipschitz constant of $f^L_k$, $N\geq0$ be the number of iterations, $t_k=t\geq1$ be the number of samples, $\mu>0$ be the  smoothing parameter in \eqref{eq:gsmooth} and $\mathcal{U}_k = [u_0,u_1,\cdots,u_k]$, 
     $k\in [N]$.
     Moreover, let $\{x_k\}_{k\geq0}$ be the sequences generated by Algorithm~\ref{alg:ZOGD} with \eqref{eq:gon} and \eqref{eq:newup} as the update step in Lines \ref{line:6} and \ref{line:5}.  
     Then, for any iteration $N$, with \(h_{k}= h = \frac{(n+3\sqrt{n}P_N^\ast)^{\frac{1}{2}}} {(N+1)^{\frac{1}{2}}L_0(n+4)}\) and \(\mu\leq\frac{1}{L_0n^{\frac{1}{2}}(N+1)^{\frac{1}{2}}}\),
     we have
	\begin{align}
    \begin{split}\label{eq:eqond}
	\sum_{k=0}^{N}&E_{\mathcal{U}_k}[f_k^L(x_k)]-\sum_{k=0}^{N}\min_{x\in\mathcal{K}}f_k^L(x) \\
    & \leq(N+1)^{\frac{1}{2}}\left(1+ (n+4)\bar{L}_0(n+3\sqrt{n}P_N^\ast)^{\frac{1}{2}}\right), 
    \end{split}
	\end{align}
    where $P_N^\ast=P_N(x^\ast_0,\dots,x^\ast_N)$. Moreover, it holds that
    \begin{align}\label{eq:eqond2}
		E_{\tau,\mathcal{U}_{N-1}}[\text{Regret}_N]= O(n^{\frac{3}{2}}\sqrt{NP^\ast_N}),
    \end{align}
    where $\tau$ is the threshold defined in Algorithm~\ref{alg:ZOGD}
\end{theorem}
\begin{proof}
    Let $e_k = \|x_k-x_k^\ast\|.$ Then, we have
    \begin{align*}
    e_{k+1}=&\|x_{k+1}-x_{k+1}^\ast + x_k^\ast - x_k^\ast \|\\
    \leq& \|x_{k+1}-x_k^\ast\| + \|x_{k+1}^\ast-x_{k}^\ast\|    
    \end{align*} 
    and with
    $\|x_{k+1}-x_k^\ast\|\leq\sqrt{n},$ 
    it holds that
    \begin{align}\label{eq:th_ond1}
        e_{k+1}^2 \leq \|x_{k+1}-x_k^\ast\|^2\!+\!\|x_{k+1}^\ast-x_{k}^\ast\|^2\!+\! 2\sqrt{n}\|x_{k+1}^\ast-x_{k}^\ast\|.
    \end{align}
    
    Using the non-expansiveness property of projections on
    convex sets, we have
    \begin{align*}
        \|x_{k+1}-x_k^\ast\|&^2\!\leq\!\|\bar{x}_{k+1}\!-\!x_k^\ast\|^2\leq\|x_k\!-\!x_k^\ast-hg_{\mu,k^+}(x_k)\|^2\notag\\
        =&e_k^2\!-\!2h\langle g_{\mu,k^+}(x_k),x_k\!-\!x_k^\ast\rangle\!+\!h^2\|g_{\mu,k^+}(x_k)\|^2.
    \end{align*}
    Then, considering \eqref{eq:th_ond1}, the estimate 
    \begin{align}\label{eq:th_ond3}
        e_{k+1}^2\leq&e_k^2-2h\langle g_{\mu,k^+}(x_k),x_k-x^\ast\rangle+h^2\|g_{\mu,k^+}(x_k)\|^2\notag\\&+\|x_{k+1}^\ast-x_{k}^\ast\|^2+ 2\sqrt{n}\|x_{k+1}^\ast-x_{k}^\ast\|
    \end{align}
    is obtained.
    Taking expectation of 
    \eqref{eq:th_ond3} with respect to $u_k$ and considering using similar steps as in
    the proof of Theorem~\ref{th:off}, we have
    \begin{align*}
        E_{u_k}[e_{k+1}^2]\leq&e_k^2-2h\langle \nabla f^L_{\mu,k}(x_k),x_k-x_k^\ast\rangle+h^2L_0^2(n+4)^2\\&+\|x_{k+1}^\ast-x_{k}^\ast\|^2+ 2\sqrt{n}\|x_{k+1}^\ast-x_{k}^\ast\|
    \end{align*}
With
the convexity of $f^L_k$, which implies that
$f_{\mu,k}^L$ is also convex,
and \cite[Thm~1]{nesterov_random_2017}, we have
    \begin{align}\label{eq:eq:th_ond4}
        &f_k^L(x_k)-f_k^L(x_k^\ast)\leq\frac{e_k^2-E_{u_k}[e_{k+1}^2]}{2h}+\mu L_0n^{\frac{1}{2}} \\
        &+\frac{hL_0^2(n+4)^2}{2}+\frac{\sqrt{n}}{2h}\|x_{k+1}^\ast-x_{k}^\ast\|+ \frac{\sqrt{n}}{h}\|x_{k+1}^\ast-x_{k}^\ast\|. \notag
    \end{align}
    Considering the fact that $\|x_0-x_0^\ast\|^2\leq n,$ taking the expectation with respect to $\mathcal{U}_{k-1}$ from \eqref{eq:eq:th_ond4}, summing from $k=0$ to $k=N,$ and substituting the given values for $\mu$ and $h,$  we arrive at \eqref{eq:eqond}. 

    From
    Lemma~\ref{lm:lmlova} we know that $f_k(S_k^\ast) = f_k^L(x_k^\ast)$ and where $S_k^\ast$ is the minimiser of $f_k(S).$
    Considering the projection step in Algorithm~\ref{alg:ZOGD}, we know that $x_k\in[0,1]^n$ for all $k.$ Thus, from
    Remark~\ref{rem:lovae}, we have that $f_k^L(x_k) = E_\tau[f_k(S_k)].$ Now, recall
    that $\tau$ and $u_k$ are independent random variables.
    Combining
    these facts and taking the expectation from $f_k(S_k)-f_k(S^\ast)$ with respect to $\tau,$ taking the expectation with respect to $\mathcal{U}_{k-1}$, summing from $k=0$ to $k=N,$ and considering \eqref{eq:eqond}, we arrive at \eqref{eq:eqond2}.
\end{proof}
In the next section, we will demonstrate the theoretical findings through a numerical example.

\section{Numerical Example}\label{sec:exm}
In this section, we will illustrate the theoretical results given in Section~\ref{sec:main} through three numerical examples{\footnote{The Python code is publicly available at \url{https://github.com/amirali78frz/Minimisation_projects.git}.}.
\subsection{Semi-Supervised Clustering}\label{sec:ssc}
In this section, we solve the problem of semi-supervised clustering on a two-moon dataset. We consider the cost function defined in \cite[Sec. 6.5]{bach2013learning}. Given $p\in \N$ data points in a certain set $V,$ we intend 
to select set $A\subset V$, which minimises the
cost function.
\begin{align*}
    Cost(A) = I(f_A,f_{V\setminus A}) \!-\! \sum_{k\in A}\log\; \eta_k \!-\! \sum_{k\in V\setminus A} \log (1-\eta_k).
\end{align*}
Here, 
$f_A$ and $f_{V\setminus A}$ are two Gaussian processes with zero mean and covariance matrices $K_{AA}$ and $K_{{V\setminus A} {V\setminus A}},$
\begin{align*}
 I(f_A,f_{V\setminus A})&=\tfrac{1}{2}\big(\log\!\det(K_{AA})+\log\!\det(K_{{V\setminus A} {V\setminus A}}) \\
 & \qquad \qquad \qquad -\log\!\det(K_{VV})\big)
\end{align*}
is the mutual information between the Gaussian processes, and $\eta_k$ is the probability of each element to be in set $A.$ 
For more information on this problem, we refer the reader to
\cite[Sec. 6.5]{bach2013learning}. We consider each data point $x_i\in\mathbb{R}^2$ and sample each data point from the “two-moons” dataset \cite{bach2013learning}. We consider 50 points and 8 randomly chosen labelled points, where we impose $\eta_k\in\{0,1\}$ for them and $\eta_k=\tfrac{1}{2}$ for the rest. We consider the Gaussian radial basis function
kernel $k(x,y) = \exp(\frac{-\|x-y\|^2}{2\sigma^2})$ with $\sigma^2=0.05$ to obtain covariance matrices. We use Algorithm~\ref{alg:ZOGD} to find  a 
minimiser of the cost function and we define the parameters 
$\mu=10^{-5},$ $h=10^{-4},$ and $N=4000$ for the simulation. Moreover, for comparison, we implement
the subgradient method described in \cite[Sec. 10.8]{bach2013learning}. In Figure~\ref{fig:fval}, we can see the Lovász extension value over the  
sequence generated  by the subgradient method and Algorithm~\ref{alg:ZOGD}. We can see that the ZO method has a similar performance as 
the subgradient method in this scenario. Although we should note that the subgradient method is at least twice as fast as
the ZO method. This is due to the fact that the calculation of the calculate the 
subgradient of the Lovász Extension requires
$n$ function queries and for the calculation of the
random oracle \eqref{eq:eq23}, we need $2n$ function queries. 
Despite this drawback in terms of the function queries, in
Figure~\ref{fig:cluster}, we
see that Algorithm \ref{alg:ZOGD} successfully recovers the clustering from a random initial condition.

\begin{figure}
    \centering
    \includegraphics[width=1\linewidth]{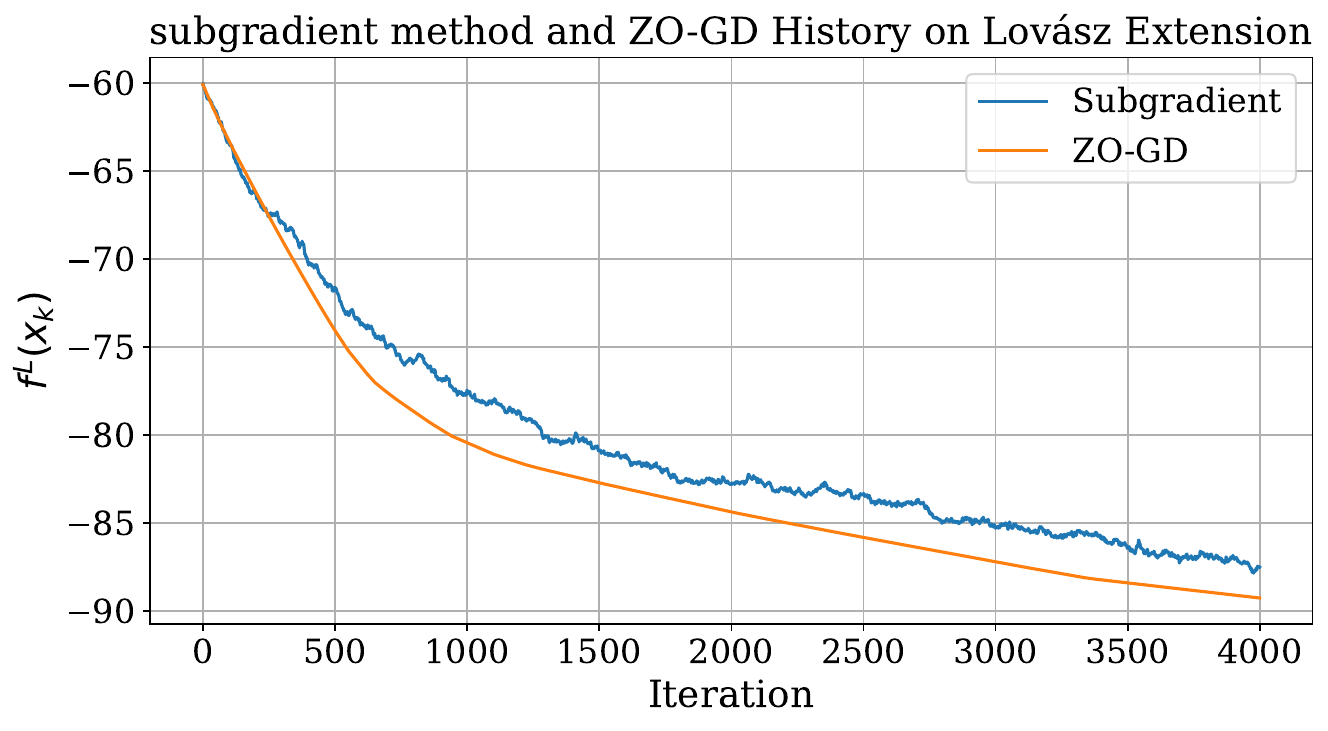}
    \caption{Lovász Extension History Over the Sequence Generated by Algorithm~\ref{alg:ZOGD} and the Subgradient Method for Offline Clustering}
    \label{fig:fval}
\end{figure}
\begin{figure}
    \centering
    \includegraphics[width=1\linewidth]{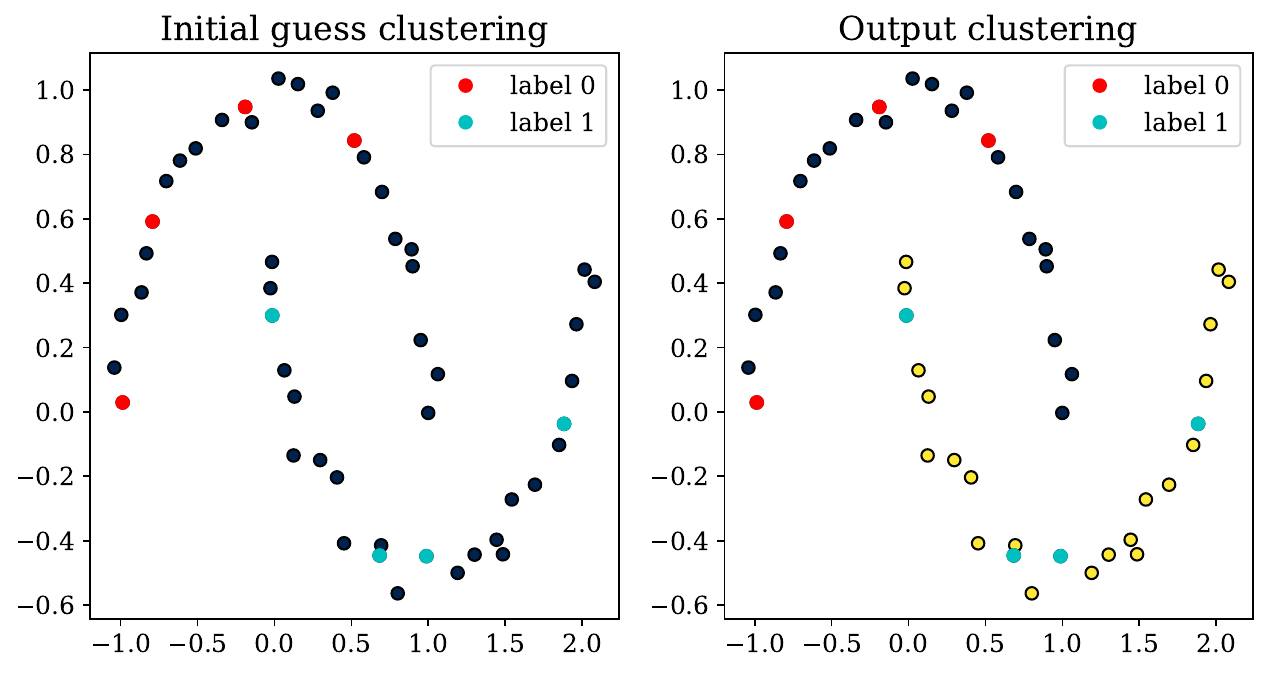}
    \caption{Offline Clustering of the Two Moon Dataset Using the Sequence Generated by Algorithm~\ref{alg:ZOGD}}.
    \label{fig:cluster}
\end{figure}

\subsection{Approximating Lovász Extension}
As discussed in the previous section, if we calculate the Lovász extension using \eqref{eq:extlov}, theoretically it can be seen that the subgradient method is at least twice as fast as
Algorithm~\ref{alg:ZOGD}. In this section we
investigate how to
calculate the Lovász extension faster, i.e., in less than 
$2n$ queries of the submodular function. Here, we consider the cost function $\text{cost}(A)=\log\!\det(K_{AA}),$ where $K_{AA}$ is defined in Section~\ref{sec:ssc}. We implement Algorithm~\ref{alg:ZOGD} and the subgradient method, as a comparison. 
We implement both algorithms using the exact and low-rank approximations of the matrix $K_{VV}$. We calculate the lower-rank approximation using the Nyström method \cite{williams2000using}.  Moreover, to evaluate the Lovász extension in Algorithm~\ref{alg:ZOGD}, we explore
two different methods based on the stochastic Lovász extension and a surrogate Lovász extension obtained through
the Taylor series approximation of $\log\!\det(K_{AA}).$
To reduce the computational cost of the Lovász extension, we adopt a stochastic sampling strategy. 
Instead of evaluating the function on all $n$ components, we randomly select a subset of components of $x$ 
according to a fixed sampling ratio $\rho \in (0,1)$. 
For a random permutation $\pi$ of $[n]$ 
and sorted weights $x_{\pi(1)} \geq \cdots \geq x_{\pi(n)}$, 
the Lovász extension is approximated by
\[
\tilde{f}^L_{\text{stoch}}(x) 
= \sum_{i \in \mathcal{I}_\rho} 
\bigl( x_{\pi(i)} - x_{\pi(i+1)} \bigr) 
\, f\bigl( S_i \bigr),
\]
where $\mathcal{I}_\rho \subset [n]$
is the sampled index set 
and $S_i = \{\pi(1),\dots,\pi(i)\}$. 
This provides a noisy but cheaper estimate of the Lovász extension.
In the Taylor surrogate method, we construct a smooth surrogate of $\log\det(K_{AA})$ using a second-order Taylor expansion 
around the identity. Given $w \in [0,1]^n$, we define the weighted kernel
\(
K_x = \operatorname{diag}(\sqrt{x}) \, K \, \operatorname{diag}(\sqrt{x}),
\)
and set $E = K_x - I$. Using the expansion
\[
\log\det(I+E) \approx \mathrm{tr}(E) - \tfrac{1}{2}\mathrm{tr}(E^2),
\]
the surrogate objective is
\[
\tilde{f}^L_{\text{Taylor}}(x) 
= \mathrm{tr}(K_x - I) 
- \tfrac{1}{2}\mathrm{tr}\bigl((K_x - I)^2\bigr),
\]
which avoids combinatorial evaluation and provides a differentiable approximation.
We implement
the mentioned methods with $N=1800$ and step size $h=10^{-4}.$ The Lovász extension value versus number of iterations and wall-clock time is given in Figures \ref{fig:comp1} and \ref{fig:comp2}. As can be seen, the performance of all the methods is similar in terms of finding the minimum value. Time-wise, we can see that using Taylor's surrogate method, the ZO method can achieve a similar runtime as the subgradient methods. 
\begin{figure}
    \centering
    \includegraphics[width=1\linewidth]{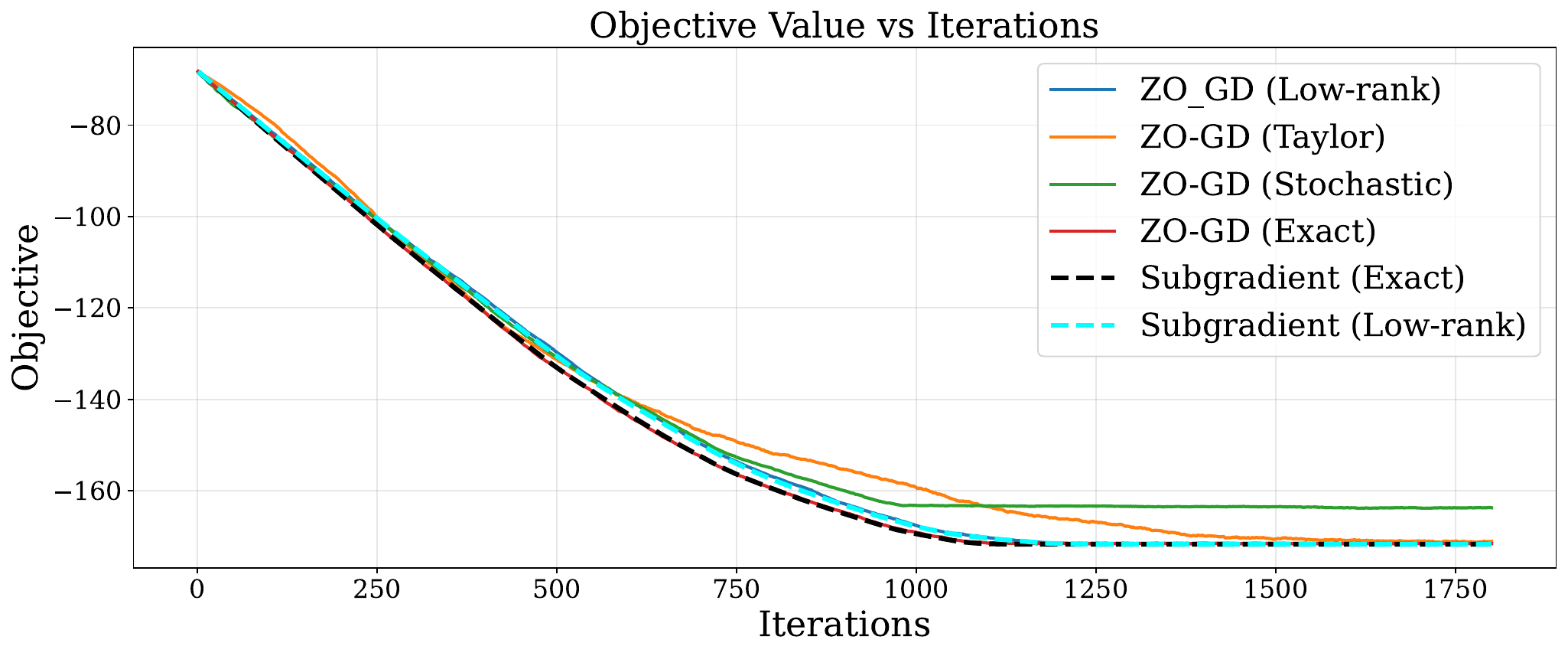}
    \caption{Logdet Cost Function's Lovász Extension History Over Iterations}
    \label{fig:comp1}
\end{figure}
\begin{figure}
    \centering
    \includegraphics[width=1\linewidth]{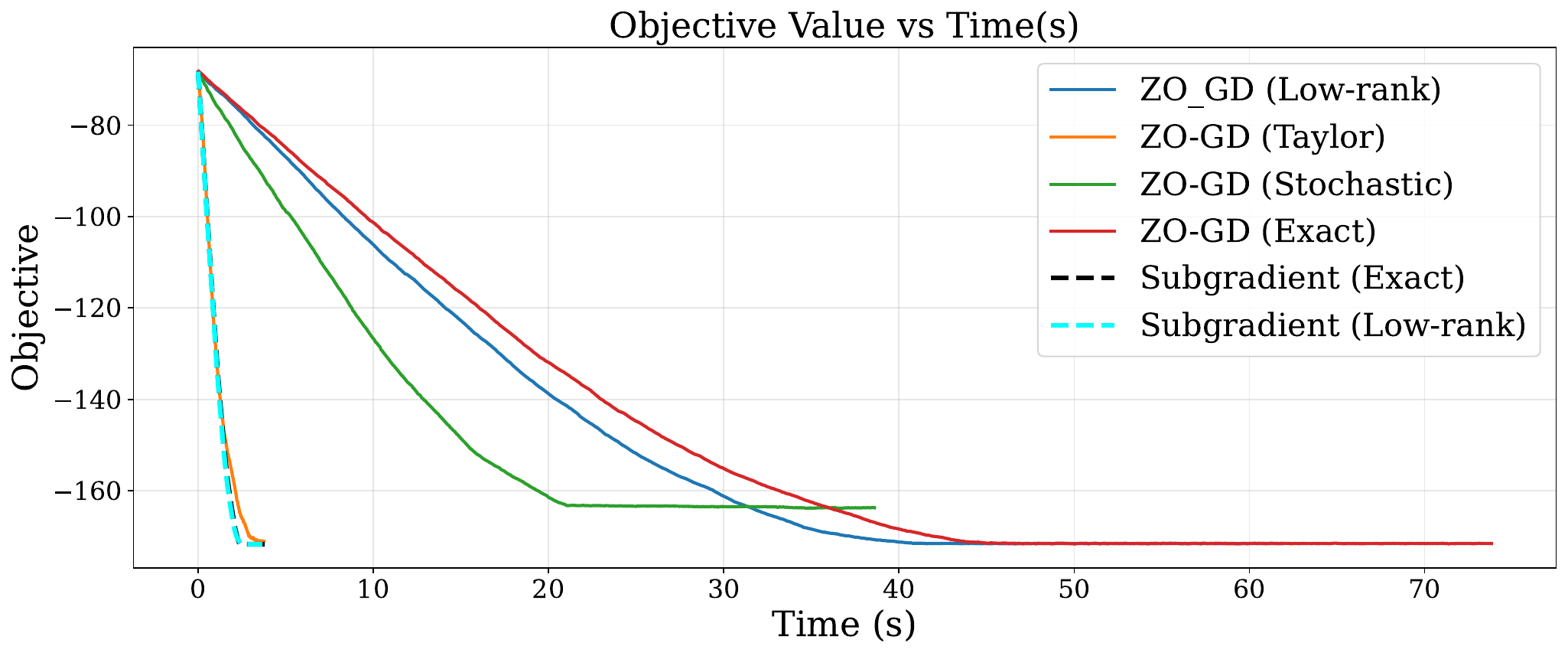}
    \caption{Logdet Cost Function's Lovász Extension History Over Time}
    \label{fig:comp2}
\end{figure}

As it can been in Section~\ref{sec:ssc}, the cost function of the semi-supervised clustering example is a combination of $\log\!\det$ functions. Thus, we investigate the above approximations in solving the semi-supervised clustering problem. In Figure~\ref{fig:cluster2}, we can see that from a random initial condition, we could successfully cluster the data using the subgradient method and using Algorithm~\ref{alg:ZOGD} with exact, Taylor surrogate, and low rank queries of the Lovász extension. We should note that the subgradient method using lower rank approximation could not recover the true clusters for this particular example. Moreover, in Figure~\ref{fig:ffval2}, we can see the Lovász extension value over the generated sequence by the mentioned methods versus.
It can be seen that using the Taylor surrogate for querying the Lovász extension, is
$90\%$ faster compared to implementation of the exact Lovász extension.

\begin{figure}
    \centering
    \includegraphics[width=0.9\linewidth]{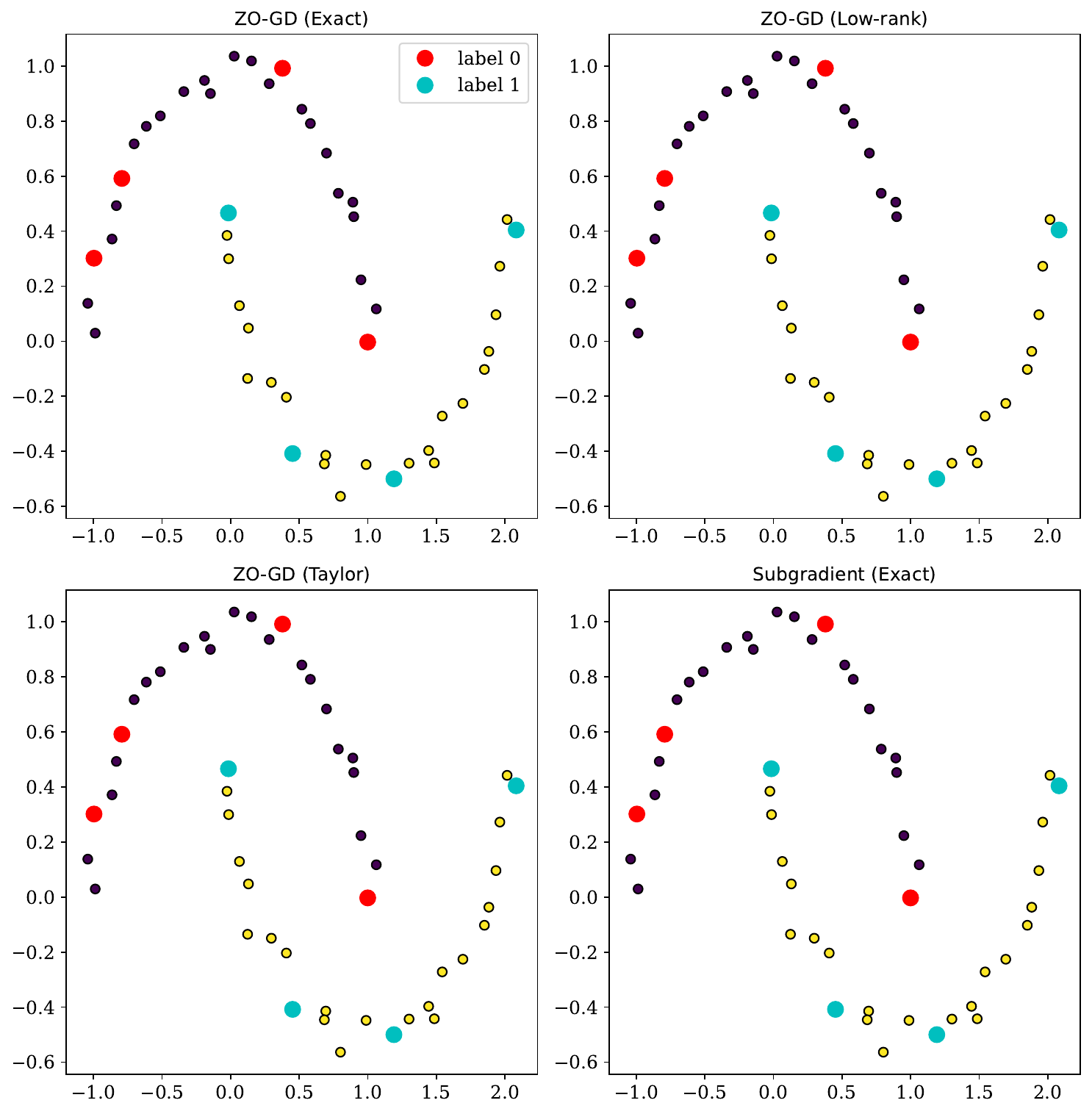}
    \caption{Offline Clustering of the Two Moon Dataset Using the Subgradient Method and Algorithm~\ref{alg:ZOGD} with Exact and Approximate Lovász Extension Queries}
    \label{fig:cluster2}
\end{figure}
\begin{figure}
    \centering
    \includegraphics[width=1\linewidth]{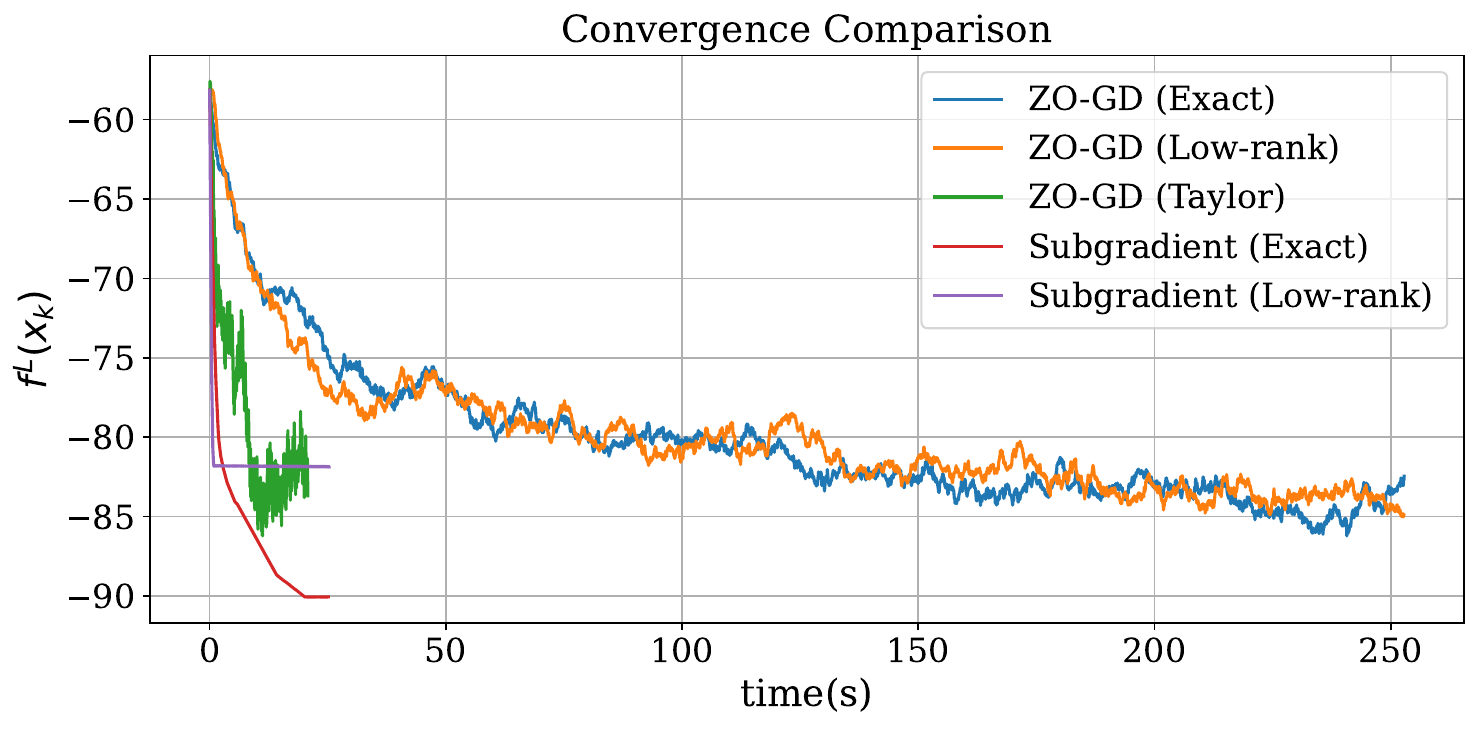}
    \caption{Lovász Extension History Over the Sequence Generated by Subgradient Method and Algorithm~\ref{alg:ZOGD} with Exact and Approximate Queries for Offline Clustering}
    \label{fig:ffval2}
\end{figure}

\subsection{Online Semi-Supervised Clustering}
In this section, we solve the online semi-supervised clustering problem where each of the true clusters is moving with a different trajectory in each iteration. Thus, $K_{VV}$, defined in Section~\ref{sec:ssc}, changes constantly, leading to an online problem.
In this section, we solve the problem using the subgradient method and Algorithm~\ref{alg:ZOGD} with the Taylor surrogate and exact queries of the Lovász extension. The Lovász extension value versus iteration and time is given in Figures \ref{fig:comp3} and \ref{fig:comp4}. Moreover, the final clustering is given in Figure~\ref{fig:cluster3} and the clustering and trajectory of the points over iterationss are given in Figure~\ref{fig:pos}. As it can be seen, all three methods have successfully completed the online clustering task while using the Taylor surrogate instead of the exact Lovász extension can make Algorithm~\ref{alg:ZOGD} up to $95\%$ faster.

\begin{figure}
    \centering
    \includegraphics[width=1\linewidth]{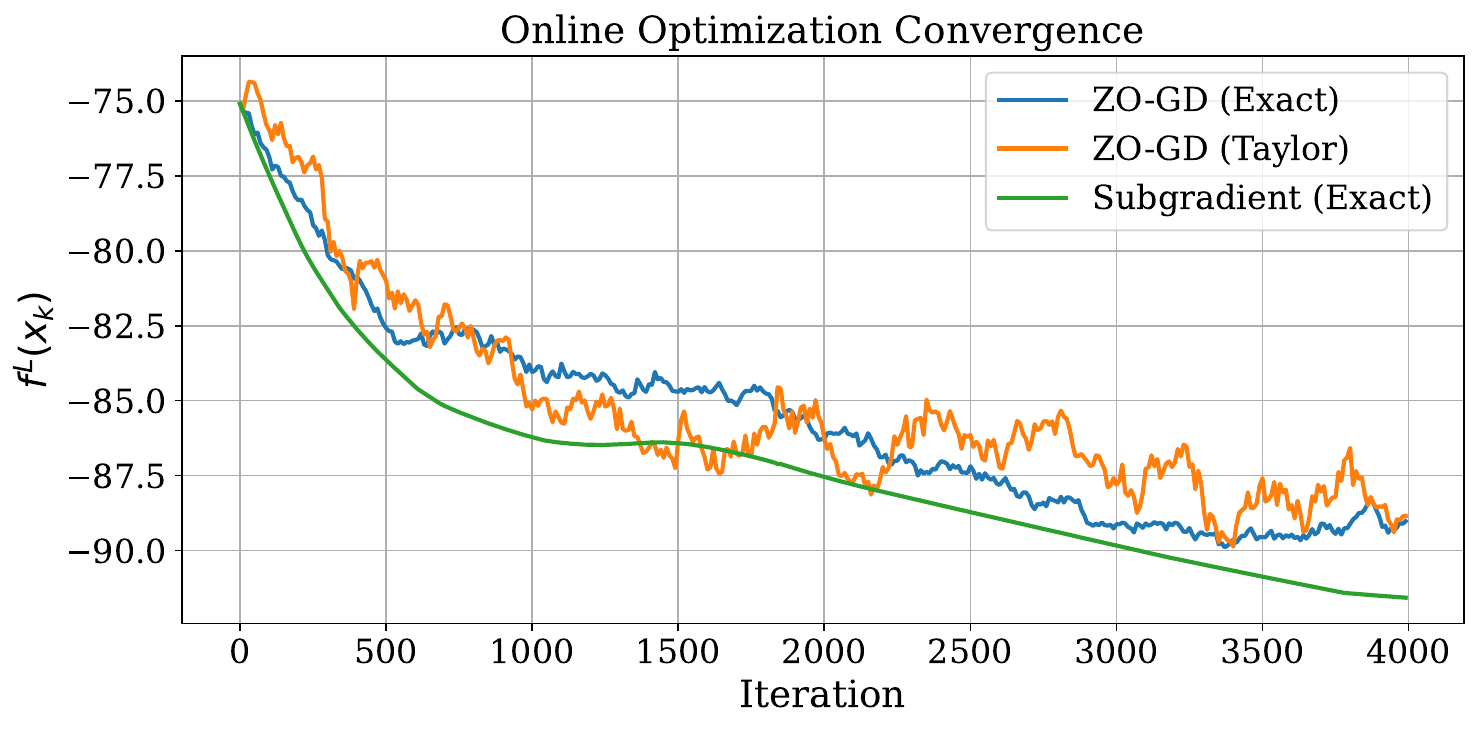}
    \caption{Lovász Extension History Over iterations for Online Clustering}
    \label{fig:comp3}
\end{figure}
\begin{figure}
    \centering
    \includegraphics[width=1\linewidth]{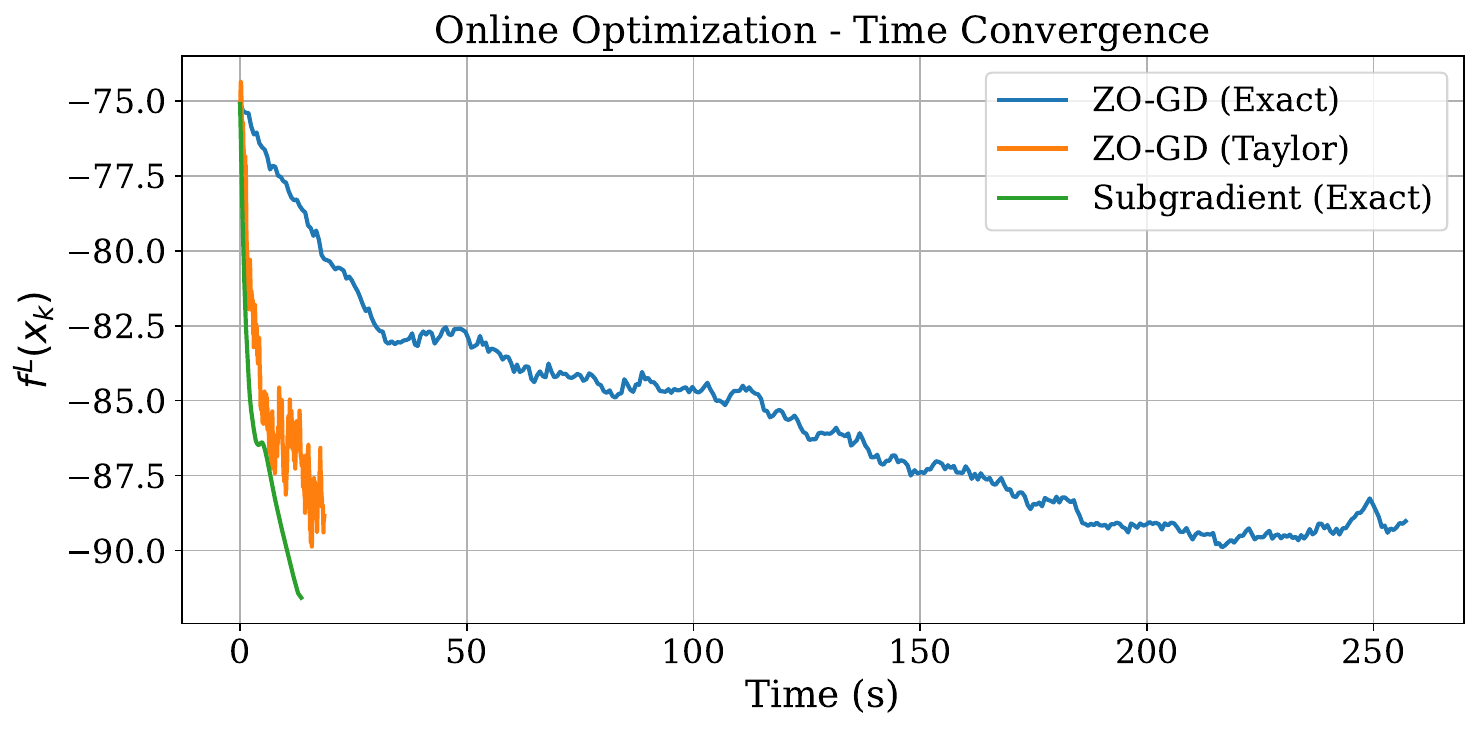}
    \caption{Lovász Extension History Over Time for Online Clustering}
    \label{fig:comp4}
\end{figure}
\begin{figure}
    \centering
    \includegraphics[width=0.9\linewidth]{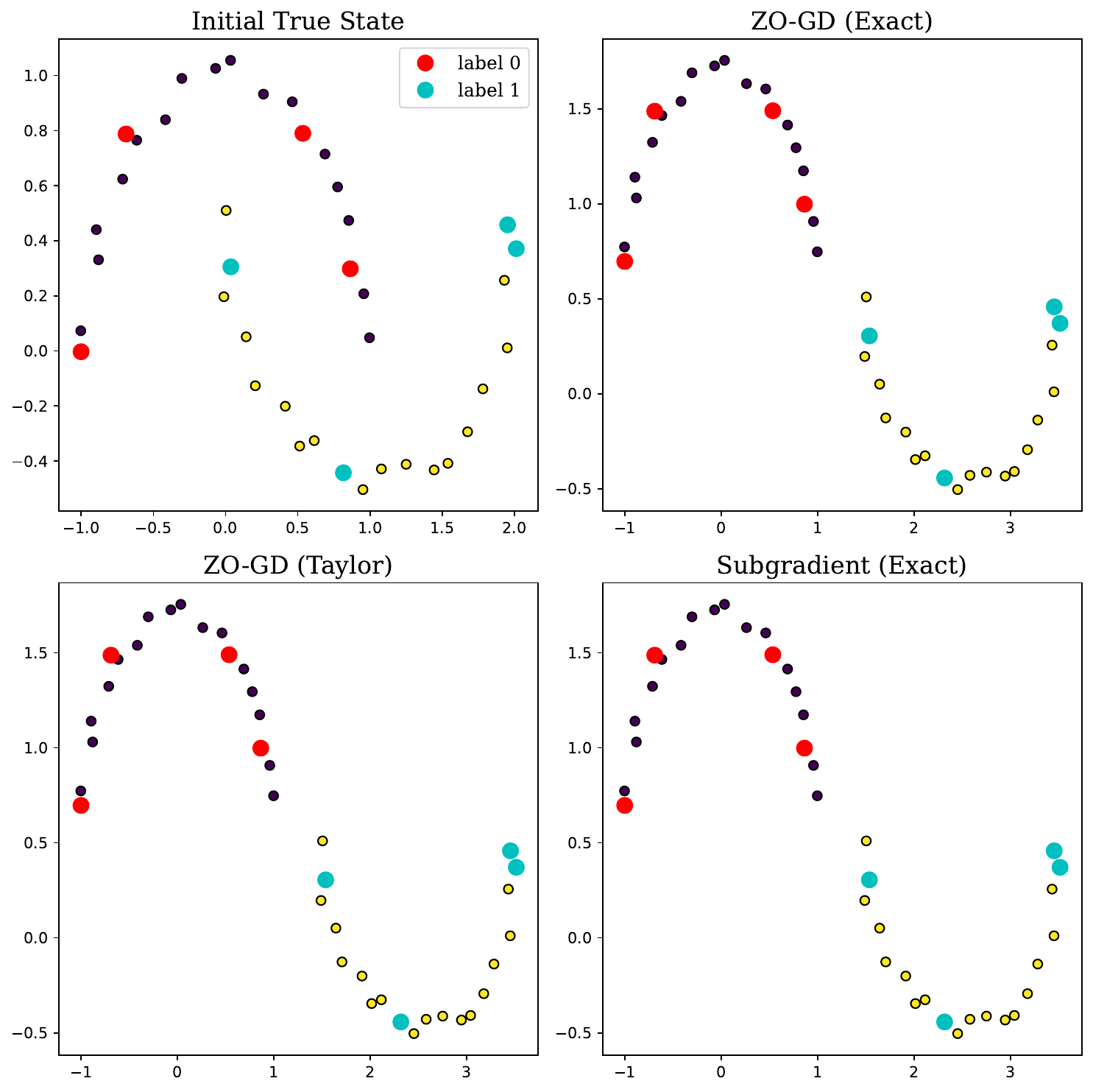}
    \caption{Online Clustering of the Two Moon Dataset Using the Subgradient Method and Algorithm~\ref{alg:ZOGD} with Exact and Approximate Lovász Extension Queries}
    \label{fig:cluster3}
\end{figure}
\begin{figure}
    \centering
    \includegraphics[width=0.9\linewidth]{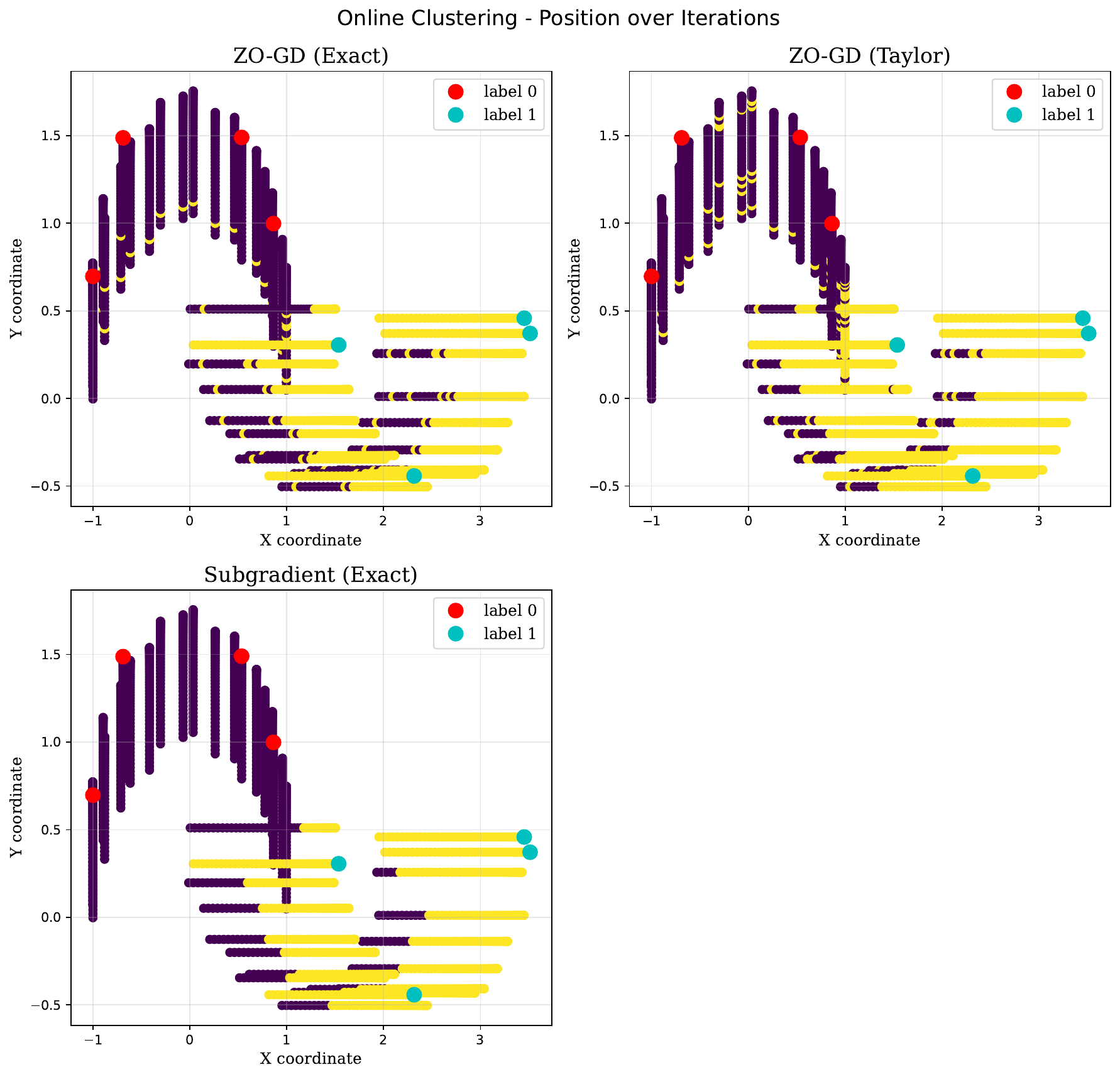}
    \caption{Positions and Labels Over Iterations for Online Clustering}
    \label{fig:pos}
\end{figure}
\section{Conclusions and Future Work}\label{sec:conc}
In this study, we considered the submodular function minimisation problem in both offline and online settings. We considered a Gaussian smoothing ZO method and investigated its performance in solving the SFM problem.
First, we considered the offline SFM and showed that the framework, in the expectation sense, converges to an $\epsilon$-approximate solution in $O(n^2\epsilon^{-2})$ function calls. Then, we considered the online SFM under both static and dynamic regrets. We showed that, in the expectation sense, the algorithm is Hannan-consistent in terms of static regret, and achieves a dynamic regret of $O(\sqrt{NP_N^\ast})$. Numerical examples demonstrating the theoretical results were presented in the end. A possible future research direction is to investigate the use of different surrogates of Lovász extension in the implementation of ZO algorithms. Another interesting possible future direction is to study the performance of ZO methods in solving minimax problems by leveraging the submodularity structure.

\addtolength{\textheight}{-12cm}   




\bibliographystyle{ieeetr}
\bibliography{ECC26}

\end{document}